\documentclass[runningheads]{llncs}
\usepackage{amssymb}
\usepackage{graphicx}
\def\XX{\mathbf{X}}
\def\vecx{\vec{x}}
\def\num{\mathtt{num}}
\def\len{\mathtt{len}}
\def\bmatrix{{{B-matrix}}}
\newcommand\Lf{\mathcal{L}_{\vec{f}}}
\newcommand\grad[1]{{\nabla{#1}}}

\begin{document}
\title{Linear Invariants for Linear Systems}
%
%
\author{Ashish Tiwari\inst{1}\orcidID{0000-0002-5153-2686}}
%
\authorrunning{A. Tiwari}
%
\institute{Microsoft, Redmond WA 98052, USA 
  \email{ashish.tiwari@microsoft.com}\\
  \url{http://www.csl.sri.com/users/tiwari}}
%
\maketitle              
\begin{abstract}
A central question in verification is characterizing when a system has 
invariants of a certain form, and then synthesizing them.
We say a system has a $k$ linear invariant, $k$-LI in short, if it has a 
conjunction of $k$ linear (non-strict) inequalities -- equivalently, an 
intersection of $k$ (closed) half spaces -- as an invariant.
We present a sufficient condition -- solely in terms of eigenvalues of
the $A$-matrix -- for an $n$-dimensional linear dynamical system to have a $k$-LI.
Our proof of sufficiency is constructive, and we get a procedure that 
computes a $k$-LI if the condition holds.
We also present a necessary condition, together with many example linear
systems where either the sufficient condition, or the necessary is tight,
and which show that the gap between the conditions is not easy to overcome.
In practice, the gap implies that using our procedure, we synthesize
$k$-LI for a larger value of $k$ than what might be necessary. 
Our result enables analysis of continuous and hybrid systems with linear
dynamics in their modes solely using reasoning in the theory of linear
arithmetic (polygons), without needing reasoning over nonlinear
arithmetic (ellipsoids). 

\end{abstract}
\section{Introduction}
\label{sec:intro}

Linear systems are extensively studied because they serve as a 
good modeling formalism. 
Even systems that have nonlinear dynamics often exhibit nice linear
behavior in certain regions of the state space, and can be modeled
using (piecewise) linear systems or hybrid systems with linear continuous
dynamics.  Furthermore, 
linear systems are easier to analyze and can be used to build analyzers for
nonlinear, piecewise linear, and hybrid systems.  

A linear system is simply a continuous-time and continuous-space dynamical
system whose state space is the $n$-dimensional reals, and whose dynamics
is specified by an ODE of the form $d\vec{x}/dt = A\vec{x}$~\cite{Kailath80}.
A classical result in control says that this system is stable (around 
the origin) if the real parts of all eigenvalues of $A$ are negative.
Moreover, for such stable $A$, there exist quadratic Lyapunov functions:
functions that are decreasing along the system trajectories.
A Lyapunov function $f$ gives rise to invariants: if a system starts 
in the region $f(\vec{x}) < 1$, then it will continue to stay inside
that region.  Such Lyapunov (and Lyapunov-like) functions have been
used extensively in verification of linear, nonlinear, and hybrid 
systems~\cite{DBLP:conf/hybrid/DuggiralaM12,boydL,emsoft13a}.

There are, however, a few undesirable features if we just rely on ellopsoidal
invariants for stable systems. First, we also want to generate invariants for unstable
systems. Second, 
while we have made significant progress in 
reasoning with nonlinear real arithmetic~\cite{Dejan12}, it is still much more scalable
and desirable to have linear (or piecewise linear) functions $f$ defining the invariant.
We are interested here in invariants of the form $f(\vec{x}) < 1$ 
that can be expressed solely using linear expressions.
In other words, we are interested in {\em{conjunction of linear inequalities}}
as invariants, and want to 
know how to compute good quality invariants of this form.  
We study this question for the class of linear systems in this paper.

We use the term {\em{$k$ linear invariant}}, or $k$-LI in short,
to denote an invariant that can be represented as a conjunction of $k$ linear
inequalities.
The parameter $k$ here provides a good tradeoff between strength of the invariant
and the cost of reasoning with it:
(1) a small value of $k$ makes the task of reasoning with these invariants easier, but
it also restricts the strength of the invariant, whereas
(2) a large value of $k$ can yield potentially strong invariants, but it also makes the reasoning
task more complex. 
Ideally, we want $k$ to be large enough to allow {\em{polytopes}} as candidate invariants,
but not much larger than needed.
Polyhedral Lyapunov functions~\cite{BlanchiniMiani},
for example, define such polytopes with $2s$ faces, where $s \geq n$.
A special case is when the polytope in $n$-dimensional space
has $2n$ faces, such as a zonotope (rotated box).  
As we will show, for stable systems, we generate $2s$-LI where $s \geq n$, whereas for
unstable systems, we can still generate $2s$-LI where $s$ may be less-than $n$.

Our main result is a sufficient condition that guarantees existence of 
$k$ linear invariants ($k$-LI) for linear systems. 
A corollary of the main result is a sufficient condition for existence of 
polyhedral Lyapunov functions for linear systems. 
Our proof is constructive:
we actually show how to generate the $k$-LI if our
condition holds.  We then go on to answer the question about whether the
condition is also necessary.  
We present examples to show the condition is not necessary, and we present
a slightly weaker necessary condition. It is an interesting open question
if we can close the gap between the sufficient condition and the 
necessary condition.
Intuitively, the sufficient condition gives us the number $k$ of half-spaces
needed to guarantee existence of a bounded $k$-LI, while the necessary condition gives 
the number $k$ of half-spaces so that there is guaranteed no $k'$-LI for $k' < k$.

Our results on necessary and sufficient conditions for existence of linear
invariants for linear systems are foundational:
the real part of eigenvalues of the $A$-matrix being negative characterizes
existence of {\em{ellipsoidal}} invariants, and here we present similar conditions for
existence of {\em{linear}} invariants.
There is a long line of work on trying to characterize linear invariants starting with
Bitsoris and Kiendl~\cite{Bitsoris91,Kiendl92:TAC}.
Our results nontrivially extend these results, and will be of interest in 
(a) verification of linear and hybrid systems, as verification techniques can be 
strengthened by computing strong linear invariants for modes with linear dynamics, 
and
(b) applications of infinity-norm, or polyhedral, Lyapunov functions (PLFs), as 
our results provide a simple procedure for generating PLFs for linear systems.


%
%
\section{Preliminaries}
\label{sec:prelims}

\subsection{Polyhedron}
We use the term {\em{polyhedron}} to mean any subset of $\Re^n$ defined by
intersection of half-spaces, and the term {\em{polytope}} to mean a polyhedron
that is also bounded.
A polyhedron that strictly contains the origin
can be written as a conjunction of linear inequalities
\begin{eqnarray*}
a_{11}x_1 + \cdots + a_{1n}x_n & \leq & 1
\\ \vdots \\
a_{s1}x_1 + \cdots + a_{sn}x_n & \leq & 1
\end{eqnarray*}
In matrix notation, this can be written succinctly as,
$F\vec{x} \leq \vec{1}$, where $F$ is an $s\times n$ matrix that
contains $a_{ij}$ in its $(i,j)$-cell.
This can be equivalently written as
$\max(F\vec{x}) \leq 1$.

A polyhedron is {\em{negative-closed}} if,
whenever
a point $\vec{c}$ occurs in the polyhedron, then
$-\vec{c}$ also occurs in it.
A negative-closed polyhedron can be represented as
$$
\max(F\vec{x}) \leq 1, \;\;
\min(F\vec{x}) \geq -1,
$$
which can be equivalently written as,
$\max(|F\vec{x}|) \leq 1$, where $|\cdot|$ is the absolute value function.
The above constraint
can be written as,
$||F\vec{x}||_\infty \leq 1$, where
$||\cdot||_\infty$ denotes the infinity-norm of a 
vector. Note that the infinity-norm of a vector is the maximum of the
absolute value of the components of the vector\footnote{%
In general, the $l$-norm $||\vec{b}||_l$ of the vector $\vec{b}$ is
$||\vec{b}||_l = (|b_1|^l + \cdots + |b_n|^l)^{\frac{1}{l}}$
where $b_1,\ldots,b_n$ are the $n$ components of the vector $\vec{b}$ and $|b_i|$ denotes
the absolute value of $b_i$.}.

If we define a function
$f: \Re^n \mapsto \Re$ such that
$f(\vec{x}) = ||F\vec{x}||_\infty$, then the 
negative-closed polyhedron 
$||F\vec{x}||_\infty \leq 1$ can be written as
$f(\vec{x}) \leq 1$.
We are interested in finding functions $f$ that define
``invariants'' for a given linear system.

If $F$ is a $s\times n$ matrix, then in general
the region $||F\vec{x}||_\infty \leq 1$ could be unbounded.
However, if the matrix $F$ has full column rank, then 
this region is bounded, and represents a {\em{polytope}} (a closed, bounded polyhedron).
%
A specific case is when $F$ is a $n\times n$ full rank matrix. 
In this case, the set $||F\vecx||_\infty\leq 1$ is a zonotope~\cite{Girard05}.
Invariant sets are closely related to Lyapunov functions, and the special case
when $F$ is $s\times n$ with $s \geq n$ and $F$ has full column rank
has been studied under the name of ``infinity-norm Lyapunov functions''
(also called polyhedral, or piecewise-linear, Lyapunov functions)~\cite{Kiendl92:TAC}.

Note that, in logical notation, 
the set $||F\vecx||_\infty \leq c$ can also be described by the formula
$$
 \bigwedge_{i=1}^{s}  (-c \leq F_{i,*}\vecx \leq c)
$$
where $F_{i,*}$ denotes the $i$-th row of the matrix $V$. The formula above is just a conjunction of
$2s$ linear (non-strict) inequalities, and can potentially be, what we will call, a $2s$-LI.

\subsection{Invariant Sets}

Consider a (continuous-time continuous-space) dynamical system with state space
$\Re^n$.  Assume that its trajectories are given by a function
$g: \Re^+ \times \Re^n \mapsto \Re^n$,
where 
$g(t,\vec{c})$ denotes the state reached at time $t$
starting from $\vec{c}\in\Re^n$ at time $0$.
Note that $g(0,\vec{c}) = \vec{c}$.

\begin{definition}[Invariant Set]
  A set of states $f(\vec{x})\leq 1$ is an
  invariant set for (a system whose trajectories are defined by) 
  $g$ if,
  for every $\vec{c}$ s.t. $f(\vec{c})\leq 1$, 
  it is the case that 
  $f(g(t,\vec{c}))\leq 1$ for all $t \geq 0$.
\end{definition}

The notion of invariance considers states reached at all future
times $t \geq 0$ -- this has also been called {\em{positive invariance}}~\cite{Blanchini99}.
For dynamical systems where $g$ is a continuous function of the
time $t$ and the initial state $\vec{c}$, 
{\em{invariance}} as defined above is equivalent to {\em{inductive invariance}}
defined below, where we only consider states reached in some small
time interval.

\begin{definition}[Inductive Invariant Set]
  A set of states $f(\vec{x})\leq 1$ is an
  inductive invariant set for $g$ if,
  for every $\vec{c}$ s.t. $f(\vec{c})\leq 1$, there exists a 
  $\delta > 0$ s.t. for all $t \in [0,\delta)$,
    $f(g(t,\vec{c}))\leq 1$.
\end{definition}

A consequence of the equivalence between
invariance and inductive invariance
for continuous dynamics is that for checking invariance,
we only need to look at points $\vec{c}$ on the boundary, that is,
$f(\vec{c}) = 0$, and not worry about trajectories starting from
points $\vec{c}$ in the interior ($f(\vec{c}) < 1$).
This fact has been used extensively in 
safety verification of continuous and hybrid dynamical 
systems~\cite{PrajnaJadbabaie04:HSCC,Tiwari03:HSCC,Platzer08:CAV,SSM04:HSCCsmall}.

We are interested in linear dynamical systems whose dynamics is specified as
$\frac{d\vec{x}}{dt} = A\vec{x}$, where $A$ is an $n\times n$ matrix.
In this case, we know that there exists a function $g$ that
satisfies the
equations $\frac{d g(t,\vec{c})}{dt} = A g(t,\vec{c})$ 
and $g(0,\vec{c}) = \vec{c}$, and that $g$ is continuous in both its
first and second arguments.

We will use the term {\em{$2s$ linear invariant, in short $2s$-LI}}, to denote
an invariant for such a dynamical system that is
of the form $||F\vecx||_\infty\leq 1$ where $F$ is an $s\times n$ matrix.



\subsection{Invariants and Weak Lyapunov Functions}

Linear invariants for linear systems are closely related to weak polyhedral
Lyapunov functions.

Let us fix a linear system
$\frac{d\vec{x}}{dt} = A\vec{x}$, and  the function
$g$ specifying its trajectories.
A (weak) Lyapunov function is (non-increasing) decreasing along all 
possible system trajectories.

\begin{definition}[PLF]\label{def:plf}
  Let $f: \Re^n \mapsto \Re^{0+}$ be the 
    function $f(\vec{x}) := ||F\vec{x}||_\infty$ where $F$ is a full column-rank
    $s\times n$ matrix.
The function $f$
    is a {\em{polyhedral Lyapunov function (PLF)}}
for $g$
if
$f(g(t,\vec{x})) < f(\vec{x})$ for all $t > 0$ and
for all $\vec{x}$ s.t. $f(\vec{x}) > 0$.
The function $f$
is a weak polyhedral Lyapunov function (weak PLF) 
for $g$
if $f(g(t,\vec{x})) \leq f(\vec{x})$ for all $t > 0$ and 
for all $\vec{x}$ s.t. $f(\vec{x}) > 0$.
\end{definition}



Existence of weak PLF is equivalent to existence of bounded linear invariants.
\begin{proposition}\label{prop:plfisinv}
Let $f: \Re^n \mapsto \Re^{0+}$ be the
    function $f(\vec{x}) := ||F\vec{x}||_\infty$ where $F$ is a full column-rank
    $s\times n$ matrix.
The function $f$ is a weak polyhedral Lyapunov function (weak PLF) 
for $g$ if, and only if, for every constant $c > 0$, 
the set $\{\vec{x} \mid f(\vec{x}) \leq c\}$
is an (inductive) invariant set for $g$.
\end{proposition}
\begin{proof}(Sketch)
If $f$ is a weak PLF, then $f$ does not increase along
trajectory starting from any point $\vec{x}$ s.t. 
$f(\vec{x}) > 0$, and hence $f(\vec{x})\leq c$ is clearly an invariant
for every $c > 0$.
Conversely, if $f$ does increase along a trajectory starting from some
point $\vec{x}_0$, then the set
$f(\vec{x}) \leq d$ is not an invariant, where $d = f(\vec{x}_0) > 0$. \qed
\end{proof}

For linear systems, the function $g$ defining the trajectory has a 
special linearity property:
$g(t,\alpha\vec{x}+\beta\vec{y}) = \alpha g(t,\vec{x}) + \beta g(t,\vec{y})$.
Using this property, 
invariance of $f(\vec{x}) \leq c$ for all $c$ is equivalent to
invariance of $f(\vec{x}) \leq 1$.

\begin{proposition}\label{prop:invisinv}
 Let $f: \Re^n \mapsto \Re^{0+}$ be as in Proposition~\ref{prop:plfisinv}.
The set $\{\vec{x} \mid f(\vec{x}) \leq c\}$
is an (inductive) invariant set for $g$ for all $c > 0$ 
if and only if
it 
is an (inductive) invariant set for $c=1$.
\end{proposition}
\begin{proof}(Sketch)
Given a point $\vec{x}$, if
$\alpha = f(\vec{x})$, then
$f(\vec{x}/\alpha) = 1$.  Hence, if a trajectory starting
from $\vec{x}_0$ exits the set $f(\vec{x}) \leq f(\vec{x}_0)$,
then the trajectory starting from 
$\vec{x}_0/f(\vec{x}_0)$ exits the set $f(\vec{x}) \leq 1$.  \qed
\end{proof}

The two propositions above together show that the ability to
compute linear invariants will also give us the ability to obtain
weak PLFs for linear systems.
Next, we show that we can restrict our focus to
{\em{negative-closed polyhedron}} without loss
of any generality.

Recall that a negative-closed polyhedron is written as 
$||F\vec{x}||_\infty\leq 1$.
The next lemma, Lemma~\ref{lemma:negclosed}, states that if there is a
polyhedral invariant then there also exists a negative-closed polyhedral invariant.
Its proof 
relies on the interesting observation that
(1) given a polyhedron $P$, we can get a smaller negative-closed 
polyhedron $Q$ contained inside $P$, and
(2) whenever
a point $\vec{c}$ is on the boundary of $Q$, then
we can prove that trajectories are pointing inwards by
either looking at where $\vec{c}$ lies on $P$ and
using PLF property of $P$, or 
looking at $-\vec{c}$ and then using PLF property of $P$.

\begin{lemma}\label{lemma:negclosed}
  If $\max(F\vec{x}) \leq 1$ is an invariant for
  a linear system $\dot{\vec{x}} = A\vec{x}$, then
  $||F\vec{x}||_\infty \leq 1$ is also an invariant.
\end{lemma}
\begin{proof}
  Let $P$ be the polyhedron $\max(F\vec{x}) \leq 1$,
  and let $Q$ be the polyhedron
  $||F\vec{x}||_\infty \leq 1$.
  Clearly, $Q \subseteq P$.
    Let $F_{1,*},\ldots,F_{s,*}$ be the rows of $F$ (thought of as
  row vectors).

  Let $\vec{c}$ be an arbitrary point in $Q$.
  We need to prove that trajectories starting from a point in $Q$ remain
  in $Q$.
%
  Since $Q$ is negative-closed,
  $\vec{c}\in Q$ implies
  $\vec{-c}\in Q$, and hence $\vec{-c}\in P$,
  and since $P$ is an invariant, we know
    $g(t,\vec{-c}) \in P$ {for all $t \geq 0$}, which means
  \begin{eqnarray}
      F_{i,*} g(t,\vec{-c}) & \leq & 1  \qquad \forall{i},\; \forall{t \geq 0}  
    \forall{\vec{c}\in Q}
    \label{eqnmc}
  \end{eqnarray}
  Using Fact~\ref{eqnmc}, we will show that 
  $g(t,\vec{c})$ is in $Q$ for
  all $t \geq 0$.

  We prove by contradiction.
  Suppose $g(t,\vec{c})$ is not in $Q$ for all $t\geq 0$.
  Define $t_\mathit{escape} = \inf\{ t \mid g(t,\vec{c})\not\in Q\}$.
  Consider the point $\vec{d} = g(t_{\mathit{escape}},\vec{c})$.
  Point $\vec{d}$ is in $Q$, but points reached from $\vec{d}$,
  that is, $g(t_{\mathit{escape}}^+,\vec{c})$, are not in $Q$, but they are in $P$;
  hence, there exists some index $i\in\{1,\ldots,s\}$ s.t.
    $F_{i,*} \vec{d} = -1$ and $F_{i,*} g(0^+,\vec{d}) < -1$.

  Consider the point $-\vec{d}$. 
    We note that $F_{i,*}(-\vec{d}) = -(-1) = 1$ by linearity.
  From Fact~\ref{eqnmc}, we get
    $F_{i,*} g(0^+,\vec{-d}) \leq  1$. 
    By linearity, this implies $F_{i,*} g(0^+,\vec{d}) \geq  -1$. 
    This contradicts $F_{i,*} g(0^+,\vec{d}) < -1$. 
This completes the proof. \qed
    %
%
%
%
 %
\end{proof}



\subsection{Necessary and Sufficient Check for Invariance}
Before we can say anything interesting about (positive) invariant sets, we need a way to
establish when some set is a (positive) invariant set and when it is not.
The Lie derivative helps here: given the vector field $\vec{f}$ (of the dynamical system), 
and a function $g: \XX \mapsto \Re$, the Lie derivative $\Lf(g)$ is defined as
the dot-product
$\nabla{g}\cdot \vec{f}$ of the gradient of $g$ with $\vec{f}$; that is,
$\Lf(g)(\vecx) := \sum_{i=1}^{n} \frac{\partial{g(\vecx)}}{\partial{x_i}}\vec{f}_i(\vecx)$.
Note that $\Lf(g)$ is just the time derivative of $g$, and the gradient $\nabla{g}$ is the 
vector $[\partial{g}/\partial{x_1},\ldots,\partial{g}/\partial{x_n}]$
of the real-valued function $g:\Re^n\mapsto\Re$. We will assume $\grad{g}$ is a {\em{row}} vector,
and so $\nabla{g}\cdot \vec{f}$ is just matrix multiplication, which we will denote by
juxtaposition (and not use $\cdot$).

Proposition~\ref{prop:invisinv} shows that checking invariance is equivalent to
checking inductive invariance. When checking inductive invariance of the set
$f(\vec{x})\leq 1$,
we first note that if point $\vec{c}$ is strictly in the {\em{interior}},
that is, $f(\vec{c}) < 1$, 
then 
$f(g(t, \vec{c})) \leq 1$ for all
sufficiently small $t > 0$.
Hence, to check if $f(\vec{x})\leq 1$ is an inductive invariant, 
we only need to
worry about points $\vec{c}$ s.t. $f(\vec{c}) = 1$; that is,
the so-called {\em{boundary points}}, and prove that the vector field
points ``inwards'' at these points~\cite{BlanchiniMiani}.
Necessary and sufficient condition for checking if a vector field points
``inwards'' were discussed in~\cite{TT09:FSTTCS}.  A {\em{sufficient}} condition
for checking that the vector field $\vec{f}$ that maps $\vecx$ to $A\vecx$
is pointing ``inwards'' into the region $g(\vecx)\leq 1$ at the point 
$\vec{c}$, where $g(\vec{c})=1$ is that
$\Lf(g)(\vec{c}) < 0$. 
In general, a {\em{necessary}}, but not sufficient, condition is that
$\Lf(g)(\vec{c}) \leq 0$. 
However, for linear dynamics and linear invariant sets, this necessary condition is
also sufficient. This fact was implicitly stated in~\cite{TT09:FSTTCS}, and we make
it explicit in Proposition~\ref{prop:check}.
We first note that if the vector field $\vec{f}$ is given by $A\vecx$,
then $\Lf(g)(\vecx)$ is simply $\grad{g(\vecx)}A\vec{x}$. 

\begin{proposition}\label{prop:check}
A polyhedral region $F\vecx \leq c$ is positively invariant 
for a linear system $\dot{\vecx} = A\vecx$ iff
for each $i\in\{1,\ldots,n\}$, it is the case that
\begin{eqnarray}
  F\vecx \leq c \; \wedge\; 
    F_{i,*}\vec{x} = c  \;\Rightarrow\; F_{i,*} A \vec{x} \leq 0
  \label{eqn:check}
\end{eqnarray}
    where $F_{i,*}$ is the $i$-th row of $F$.
\end{proposition}
\begin{proof}
    The result in~\cite{TT09:FSTTCS} showed that the necessary check 
    $F_{i,*} A \vec{x} \leq 0$
    also sufficient if we can prove that the gradient is not zero on the boundary points.
    This is clearly the case for linear invariants. \qed
\end{proof}
The reason why $\Lf{g}(\vec{c}) \leq 0$ is not sufficient is that, in general, if the
``first-derivative'' is zero, we need to check the sign of the ``second-derivative'', and
if that is zero, then the sign of the ``third-derivative'', and so on~\cite{TT09:FSTTCS}. 
For linear dynamics
and linear invariants, it is also possible to prove by first principles that these additional
checks are implied by the necessary condition.
Note that the sufficient check based on strict inequalities on the right-hand side are often
used for the general case, for example, in definition of 
Barrier certificates~\cite{PrajnaJadbabaie04:HSCC}.

\subsection{Existence and Synthesis Problem}

We now formally state the problem we solve in this paper.

\begin{definition}[$2s$-LI Decision Problem]
Given a linear system $d \vec{x}/dt = A\vec{x}$ and a natural number $s > 0$, 
determine if there exists an $s\times n$ matrix $F$ 
such that 
$||F\vec{x}||_\infty \leq 1$ is an invariant of the linear system.
\end{definition}

\begin{remark}
The $2s$-LI decision problem insists on finding a 
negative-closed polyhedron.
There is no loss of generality here, since 
Lemma~\ref{lemma:negclosed} 
showed that if there is any polyhedral linear invariant 
for a linear system, then
there is one that is negative-closed. \qed
\end{remark}

The synthesis problem asks us to generate the invariant for the given linear system.

\begin{definition}[$2s$-LI Synthesis Problem]
Given a linear system $d \vec{x}/dt = A\vec{x}$, find a natural number $s > 0$ and
a matrix $F\in \Re^{s\times n}$ such that 
$||F\vecx||_\infty \leq 1$ is an invariant.
\end{definition}

\begin{remark}
While the formulation of the synthesis problem leaves open
the choice of $s$, ideally we want $s$ to be as close to $n$
as possible.
Note that there are systems that have an $2s$-LI with $s > n$,
but that have no $2s$-LI for $s \leq n$.
One such example is shown in Example~\ref{example-nonsquarePLF}.
\end{remark}

\begin{example}\label{example-nonsquarePLF}
  Consider the linear system 
  $$
   \frac{dx}{dt} \;\; =\;\; -x - \sqrt{3}y,
   \qquad
   \frac{dy}{dt} \;\; =\;\; \sqrt{3}x - y
  $$
  This models a spiral converging to the origin. 
  It has a $2s$-LI for $s=3$, but not for $s=2$; that is, it has a 
  hexagon-shaped invariant (with $6$ edges), but no parallelogram-shaped one
    (with $4$ edges).
  For example, consider the following hexagon
  $$
   |2y| \leq 1,
   \quad
   |y + \sqrt{3}x| \leq 1, 
   \quad
   |y - \sqrt{3}x| \leq 1
  $$
  Using Proposition~\ref{prop:check}, we can verify that the hexagon
  is indeed an inductive invariant. For example, the check
  in Equation~(\ref{eqn:check}) when instantiated for $i=3$ gives
  $$
   |2y| \leq 1 \; \wedge\;
   |y + \sqrt{3}x| \leq 1 \;\wedge\;
   y - \sqrt{3}x = 1 \;\;\Rightarrow\;\;
   2\sqrt{3}x+2y \leq 0
  $$
  which is actually a valid formula 
  (because $2\sqrt{3}x+2y = -2(y-\sqrt{3}x)+4y \leq -2+2\leq 0$).
  We state without proof that this system can not have a 
  parallelogram as an invariant, but we provide some intuition
    in the next example.  \qed
\end{example}

\begin{example}\label{ex:circle}
  Generalizing from Example~\ref{example-nonsquarePLF},
  it is instructive to consider the family of linear systems:
  $$
   \frac{dx}{dt} \;\; =\;\; -x - a y,
   \qquad
   \frac{dy}{dt} \;\; =\;\; a x - y
  $$
  where $a > 0$ is a positive real number.
  As $a$ becomes large and tends to infinity, the dynamics
  gets closer and closer to circular motion. 
  For circular motion, there is no $2s$-LI (as we would need
  a polyhedron with infinitely many edges to properly contain 
  a circle).  When $a = 1$, we need $4$ edges to construct an
  invariant, and when $a = \sqrt{3}$, we need $6$ edges (as in 
  Example~\ref{example-nonsquarePLF}).
  Intuitively, as $a$ increases, we need more and more edges 
  in the polyhedron to get an invariant.   \qed
\end{example}

If we solve the $2s$-LI decision and synthesis problems, then we also
get a solution for the corresponding problems for PLFs since the only
additional check needed there is to ensure that the matrix $F$ has
full column rank.

\section{Motivation and Related Work}

The decision and synthesis problems for $2s$-LI have a long history.
The motivation mainly comes from interest in {\em{polyhedral Lyapunov functions}}
(PLFs)~\cite{bl1,bl2}
to achieve robust control.
There is also work that shows that quadratic Lyapunov functions
are insufficient for establishing structural stability of 
basic motifs of biochemical networks, whereas polyhedral Lyapunov 
functions are enough to do so (of course, the motifs here 
have nonlinear dynamics)~\cite{blgi4,blgi,blgi2,blgi3}.  
Robust control and structural stability are both related to 
absolute stability -- where the goal 
is to prove stability for not a single system, but for a whole class of
systems obtained by varying some parameters in the system definition~\cite{Polanski2000}.
Informally, PLFs are a useful tool for performing such analysis, and there are also
results that show they are complete, whereas quadratic Lyapunov functions are not, 
in some cases~\cite{PLFUniversal,bl2}.
Our immediate motivation comes more from 
the use of linear invariants for analysis and verification.

Blanchini~\cite{BlanchiniNotes,BlanchiniMiani} mentions that {\em{checking}} if a 
given polyhedron is an invariant (or defines a PLF) is simple, but {\em{synthesis}} of 
such an invariant (or PLF) is not so simple.  There are iterative procedures that have 
been suggested for synthesis or procedures based on somehow ``guessing'' vertices of
the polygon~\cite{Lazar10:CDC,Lazar11:CDC,Lazar10:HSCC,Polanski2000}. However, we are
more interested in direct methods based on looking at eigenvectors of the $A$ matrix.
Bitsoris\cite{Bitsoris88,Bitsoris91} presented some of the earliest results in this
direction by giving sufficient conditions for existence of linear invariants,
and so did Kiendl~\cite{Kiendl92:TAC}. However, those works covered just the case 
when $|b| < |a|$ for a complex eigenvalue $a+b\iota$. There were no known results for
existence of linear invariants when $|b| > |a|$, which is a gap we fill in this paper.

Polyhedral invariants are essentially box invariants~\cite{AbateTiwari06:ADHS,AbateTiwariSastry09:Aut} 
but with a change of coordinates.  
A box invariant is an invariant of the form
$\max(|x_1|,\ldots,|x_n|) \leq 1$.  It is a $2n$-LI.
While it is easy to characterize existence
of box invariants based on simple checkable properties of the $A$-matrix, 
we are not aware of any such characterization for polyhedral invariants.
The difficulty comes from the fact that we have to also discover the 
transformation that we can apply to the system to turn it into one that has
box invariants. 





\section{Synthesizing Linear Invariants}

Based on our formulation of the decision and synthesis problem
for linear invariants, we henceforth restrict
ourselves to invariants of the form
$||F\vec{x}||_\infty \leq 1$, where $F$ is
an $(s\times n)$-matrix.

Our main result is a sufficient condition 
for the existence of $2s$-LI.  The result can be used to both determine if a 
linear system has an invariant, and also synthesize it. We formulate the result
for the case when $s = n$, but the proof will show that it can be used to synthesize
invariants with $s \neq n$. 
\begin{theorem}[Sufficient condition for existence of $2n$-LI]\label{theorem:main}
The linear system $d{x}/dt = A {x}$,
where $A$ is a $n\times n$ matrix, has a $2n$-LI 
if 
\begin{itemize}
  \item[(a)] $\lambda \leq 0$ for every real eigenvalue $\lambda$ of $A$
  \item[(b)] $\lambda < 0$ for every real eigenvalue $\lambda$ of $A$ such that
    $\mathtt{geom}(\lambda) < \mathtt{algm}(\lambda)$ 
\item[(c)]
$a + |b| \leq 0$ for every complex eigenvalue $(a+\iota b)$ of $A$, and
\item[(d)]
$a + |b| < 0$ for every complex eigenvalue $(a+\iota b)$ of $A$ 
such that $\mathtt{geom}(a+\iota b) < \mathtt{algm}(a+\iota b)$
\end{itemize}
\end{theorem}
Here $\mathtt{algm}(l)$ denotes the algebraic multiplicity of the
eigenvalue $l$ (that is, multiplicity of $l$ as a root of the characteristic
polynomial of $A$), and 
$\mathtt{geom}(l)$ denotes the geometric multiplicity of the eigenvalue
$l$ (that is, the maximum number of linearly independent eigenvectors
corresponding to $l$; equivalently, the dimension of the kernel of the
matrix $A - l I$). Note that $\mathtt{geom}(l) \leq \mathtt{algm}(l)$.

Given a rational matrix $A$, the conditions~(a)--(d) are decidable, and hence
it follows that we have a sound check for existence of $2n$-LIs.
Our proof of sufficiency is constructive -- it will explicitly generate the $2n$-LI
when the condition holds.
We will prove the above result by proving it for special cases in different lemmas
and finally we will put it all together.

%
%


First we consider the case when $A$ has only real eigenvalues
and prove sufficiency for such $A$.  

\begin{lemma}\label{lemma:jordan_real_1}
The linear system $d{x}/dt = \lambda {x}$,
where $\lambda \leq 0$, has a $2$-LI.
\end{lemma}
\begin{proof}
  The invariant is given by $|x| \leq 1$. \qed
\end{proof}

An $(A,I)$-Jordan block is a matrix with a submatrix $A$ on all its diagonal
positions and identity matrix $I$ on its top off-diagonal and $0$ elsewhere;
that is, a Jordan block is of the form 
\begin{eqnarray*}
   \left(
     \begin{array}{c@{\quad}c@{\quad}c@{\;\;}c@{\;\;}c@{\;\;}c@{\quad}c@{\quad}c}
      A & I & 0 & 0 & 0 & 0 
      \\
      0 & A & I & 0 & 0 & 0
      \\
      0 & 0 & A & I & 0 & 0
      \\*[-0.5em]
      \vdots & & & \ddots & \ddots & 
      \\
      0 & 0 & 0 & 0 & A & I
      \\
      0 & 0 & 0 & 0 & 0 & A 
    \end{array}
   \right)
\end{eqnarray*}

\begin{lemma}\label{lemma:jordan_real_n}
The linear system $d{x}/dt = J {x}$,
where $J$ is a $m\times m$  $(\lambda,1)$-Jordan block with
diagonal $\lambda < 0$, has a $2m$-LI.
\end{lemma}
\begin{proof}
  Let $x_1, \ldots, x_m$ denote the $m$ variables.
  We show that 
  $$\max(|x_m|,
    |\lambda x_{m-1}|,
    |\lambda^2 x_{m-2}|,\ldots,
    |\lambda^{m-1} x_{1}|) \leq 1$$ is a $2n$-LI.
    Let us denote the above formula by $\phi$.
  Consider the polyhedron defined by $\phi$, and 
  consider the face  $|\lambda^i| x_{m-i} = 1$ on this polyhedron.
  We want to prove that the direction of flow at points on this face
  is inwards; that is,
  $$
   \phi \;\wedge\; |\lambda^i| x_{m-i} = 1 \;\;\Rightarrow\;\;
    \frac{d|\lambda^i| x_{m-i}}{dt} \leq 0 
  $$
  For points on the face, we note that
  \begin{eqnarray*}
    \frac{d|\lambda^i| x_{m-i}}{dt} & = &  
|\lambda^i| (\lambda x_{m-i} + x_{m-i+1}) 
\\ & \leq &
      \lambda (|\lambda^i| x_{m-i}) + |\lambda|(|\lambda^{i-1}|x_{m-i+1}) 
\\ & \leq &
      \lambda (1) + |\lambda|(1) \; \leq \; \lambda + |\lambda| \; = \; 0
 \end{eqnarray*}
The above proof shows that the flow points inwards on all faces, and 
hence $\phi$ is a $2m$-LI. \qed
\end{proof}


We next move to the case when $A$ has complex eigenvalues.
First we start with a $2$-d case that has one complex eigenvalue.
\begin{lemma}\label{lemma:jordan_complex_1}
    The linear system
 $d{x}/dt = A {x}$, where
 $A = [a, b; -b, a]$ and $a+|b|\leq 0$, has a $4$-LI.
\end{lemma}
\begin{proof}
  It is easily verified that 
  the polyhedron $\max( |x_1|, |x_2|) \leq 1$ is a $4$-LI for this system. \qed
\end{proof}

In the next lemma, we consider complex eigenvalues in the diagonal 
of a Jordan block.
\begin{lemma}\label{lemma:jordan_complex_n}
The linear system $d{x}/dt = J {x}$,
where $J$ is a $2m\times 2m$ $(A,I)$-Jordan block 
with
$A = [a,b;-b,a]$ such that $a+|b| < 0$, 
and $I=[1,0;0,1]$, has a $4m$-LI.
\end{lemma}
\begin{proof}
  We can assume without loss of generality that $a < 0$ and $b \geq 0$.
  Define $\lambda = -a - b$.  By assumption, $\lambda  > 0$.
  Let $x_1, y_1, \ldots, x_m,y_m$ denote the $2m$ variables.
  Consider the polyhedron defined by the following constraints:
\[
  \begin{array}{rcl@{\qquad}rcl}
    |x_m|  & \leq & 1
    &
    |y_m|  & \leq & 1
    \\
    |\lambda x_{m-1}| & \leq & 1 
    &
    |\lambda y_{m-1}| & \leq & 1 
                   \\ & \vdots &  & & \vdots & \\
    |\lambda^{m-1} x_{1}| & \leq & 1 &
    |\lambda^{m-1} y_{1}| & \leq & 1 
  \end{array}
\]
  Let $\phi$ denote the conjunction of the above formulas.
  We need to prove that flows point inwards on the faces of the
  polyhedron $\phi$. Consider the face 
    $\lambda^i y_{m-i} = 1$ for an aribitrary $i$. 
    We need to show that
  $$
   \phi \;\wedge\; \lambda^i y_{m-i} = 1 \;\;\Rightarrow\;\;
   \frac{d (\lambda^iy_{m-i})}{dt} \leq 0
  $$
  Assuming $\phi$ and assuming $\lambda^i y_{m-i}=1$, we note that
  \begin{eqnarray*}
    \frac{d (\lambda^iy_{m-i})}{dt} & = &  
 \lambda^i (-b x_{m-i} + a y_{m-i} + y_{m-i+1}) 
 \\ & \leq &
 -b \lambda^i x_{m-i} + a + \lambda^i y_{m-i+1}
 \\ & \leq &
 -b \lambda^i x_{m-i} + a + \lambda 
 \\ & \leq &
 b + a + \lambda  \;\; \leq \;\; 0
  \end{eqnarray*}
It can be similarly verified that the flow points inwards for all other faces,
and thus the formula $\phi$ is an invariant. \qed
\end{proof}

Now, we can put all the above lemmas together to get a proof of Theorem~\ref{theorem:main}.
%
%
%
\begin{proof}[Theorem~\ref{theorem:main}]
  Transform $A$ into Jordan normal form $J$ and say $A = UJU^{-1}$.
  If we define $\vec{y} = U^{-1}\vec{x}$, then
  $\dot{\vec{y}} = U^{-1}A\vec{x} = U^{-1}AU\vec{y} = J\vec{y}$.
  Let $J_1,\ldots,J_k$ be the $k$ different Jordan blocks in $J$.
  For each Jordan block $J_i$, we generate an invariant $\phi_i$, and then 
  we get an invariant for $J$ by putting together the $\phi_i$'s.
  Specifically, \\
  (1) if $J_i$ is a $1\times 1$ matrix, then we use
  Lemma~\ref{lemma:jordan_real_1} to construct an invariant
  $\phi_i$ for $J_i$;
  \\
  (2) if $J_i$ is a $m\times m$ $(\lambda,1)$-Jordan block ($m > 1$),
  then we use
  Lemma~\ref{lemma:jordan_real_n} to construct an invariant 
  $\phi_i$ for $J_i$;
  \\
  (3) if $J_i$ is a $2\times 2$ matrix and corresponds to complex
  eigenvalue $a+\iota b$, then we use
  Lemma~\ref{lemma:jordan_complex_1} to construct an invariant 
  $\phi_i$ for $J_i$; and finally,
  \\
  (4) if $J_i$ is a $2m\times 2m$ $(A,I)$-Jordan block corresponding to 
  a complex eigenvalue $a+\iota b$ with $m > 1$, then we use
  Lemma~\ref{lemma:jordan_complex_n} to construct an invariant $\phi_i$ for $J_i$.

  Consider $\phi :=  \phi_1 \wedge \phi_2 \wedge \cdots \wedge \phi_k$.
  It is easy to see that $\phi$ is a $2n$-LI for $\dot{\vec{y}} = J\vec{y}$.

  Now, we can get an invariant for our original system by transforming
    $\phi$ back to original coordinates, namely, we transform $\phi(\vec{y})$ to
  $\phi(U^{-1}\vec{x})$ to get the required invariant. \qed
\end{proof}

The proof of Theorem~\ref{theorem:main} shows how one could synthesize 
linear invariants: every eigenvalue that satisfies one of the conditions in
Theorem~\ref{theorem:main} gives rise to a conjunct in the  invariant. We do not
need {\em{every}} eigenvalue to satisfy one of those conditions.

\begin{theorem}[Sufficient condition for existence of $2s$-LI]\label{thm:corollary}
The linear system $d{x}/dt = A {x}$,
where $A$ is a $n\times n$ matrix, has a $2s$-LI 
for a value $s$ obtained by adding $\num(\lambda)$ for each distinct eigenvalue
    $\lambda$ of $A$, where:
\begin{itemize}
  \item[(a)] $\num(\lambda) = \mathtt{algm}(\lambda)$ if $\lambda$ is real and $\lambda < 0$
  \item[(b)] $\num(\lambda) = \mathtt{geom}(\lambda)$ if $\lambda$ is real and $\lambda = 0$
  \item[(c)] $\num(\lambda) = \mathtt{algm}(\lambda)$ if $\lambda := a + \iota b$ is complex and
      $a + |b| < 0$ 
  \item[(d)] $\num(\lambda) = \mathtt{geom}(\lambda)$ if $\lambda := a + \iota b$ is complex and
      $a + |b| = 0$ 
  \item[(e)] $\num(\lambda) = 0$ otherwise
\end{itemize}
\end{theorem}
The proof of Theorem~\ref{theorem:main} also serves as a proof of Theorem~\ref{thm:corollary}.

\section{Generalized Sufficient Condition and a Necessary Condition}

The sufficient condition in Theorem~\ref{theorem:main} for existence of $2n$-LI 
intuitively appears to be necessary too. However, it is not. 

\begin{example}\label{example:main}
Consider the linear system whose $A$ matrix is
    \begin{equation}
        A = \left( \begin{array}{rrr} -0.5, & \sqrt{3}/2, & 0 
        \\ -\sqrt{3}/2, & -0.5, & 0 \\   0, & 0, & -2 \end{array} \right)
    \end{equation}
The matrix $A$ has one real eigenvalue, $-2$, and a pair of complex conjugate eigevalues,
    $-0.5 \pm \sqrt{3}/2\iota$. The complex eigenvalues do not satisfy the condition
    $a + |b| \leq 0$ because $-0.5 + \sqrt{3}/2 > 0$.
    Now, apply a ``change of coordinates'' transformation by defining new variables
    $y_1,y_2,y_3$ in terms of the old variables $x_1,x_2,x_3$ as follows:
    \begin{eqnarray*}
        y_1 & = & -\sqrt{3}x_1/6 - x_2/2 - \sqrt{3}x_3/3,\\
        y_2 & = &  \sqrt{3}x_1/6 - x_2/2 + \sqrt{3}x_3/3,\\
        y_3 & = &  \sqrt{3}x_1/3 - \sqrt{3}x_3/3
    \end{eqnarray*}
    If we transform the dynamical system into these new coordinates, the new $A$ matrix will
    be given by $FAF^{-1}$, which turns out to be
    \begin{equation}
        A' = \left( \begin{array}{rrr} -1, & 1, & 0 
        \\ 0, & -1, & 1 \\   -1, & 0, & -1 \end{array} \right)
    \end{equation}
    This new system has the invariant $\max(|y_1|,|y_2|,|y_3|) \leq 1$.
    Thus, we can get a 
    $2n$-LI on the $x_i$'s by replacing each $y_i$ in this invariant by its definition.  \qed
\end{example}

A square matrix $A$ is diagonally-dominant if, for each row, the (absolute value of the) diagonal is more than the sum 
of the (absolute values of the) other elements in the row; that is,
for all $i$, $|a_{ii}| \geq \sum_{j\neq i}|a_{ij}|$.

\begin{definition}[\bmatrix]\label{def:Bmatrix}
A diagonally-dominant matrix $A$ is a {\em{\bmatrix}} if its diagonal elements are all not positive; that is, for all $i$, $a_{ii} \leq 0$.
\end{definition}
The matrix $A'$ in Example~\ref{example:main} is a \bmatrix.

Proposition~\ref{prop:check} presented a necessary and sufficient condition for 
a linear system to have a polyhedral invariant.
It said that $\dot{\vecx} = A\vecx$ has a polyhedral invariant iff there is an
$(s\times n)$-matrix $F$ s.t. the $s$ universally quantified  formulas in Equation~\ref{eqn:check}
in Proposition~\ref{prop:check} are valid.
We can turn these validity checks of {\em{for-all}} formulas into checks for {\em{exists}} formulas
using Farkas Lemma to get the following result.

\begin{proposition}[Necessary and Sufficient Check]\label{prop:check2}
A linear system $\dot{\vecx} = A\vecx$ has a polyhedral invariant 
$||F\vecx||_\infty \leq 1$, where $F$ is an $(s\times n)$-matrix
    iff there exists an $(s\times s)$ \bmatrix\ $X$ such that $FA = XF$. 
\end{proposition}
\begin{proof}(Sketch)
    Informally, the $i$-th row of $FA$ represents the Lie derivative of 
    $F_{i,*}\vecx$, and since Equation~\ref{eqn:check} requires that the derivative
    be non-positive (or non-negative), Farkas Lemma says that we should be able to 
    write it as a linear combination of the rows of $F$ in a certain way.  The rows of 
    $X$ define this linear combination. The fact that $X$  is \bmatrix\  ensures that the 
    Lie derivative will be non-positive (or non-negative) at appropriate points. \qed
\end{proof}

One way to interpret Proposition~\ref{prop:check2} is that $F$ define a ``change of variables''--
it defines $s$ new variables -- each as a linear combination of the original $n$ variables.
On these new variables, call them $y_1,\ldots,y_s$, we want the $A$-matrix of system dynamics to 
be given by a \bmatrix. In that case, we get $\max(|y_1|,\ldots,|y_s|) \leq 1$ as an invariant.
Invariant of this form have been called box invariants~\cite{AbateTiwari06:ADHS,AbateTiwariSastry09:Aut}.
Thus, polyhedral invariants are just box invariants after a suitable transformation. 

In the special case of $2n$-LI (that is, $s = n$),
Theorem~\ref{theorem:main} attemted to characterize linear systems that have $2n$-LIs based 
on the eigenstructure of the $A$ matrix. However, it only provided a sufficient condition.
In this case, since $s = n$, if $F$ is full rank, then $FA = XF$ is the equivalent to 
$X = FAF^{-1}$, which implies that $X$ is similar to $A$.
Finding a necessary and sufficient condition would be equivalent to the problem of characterizing
when a matrix is similar to a \bmatrix\ completely in terms of its eigenstructure.  
There is no known solution to this problem.

\begin{figure}[t]
    \begin{tabular}{c@{$\quad$}l}
        \begin{minipage}[c]{0.20\textwidth}
            \centering
            \includegraphics[width=1.0\textwidth]{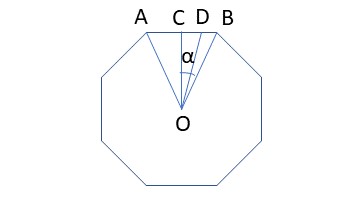}
            \caption{Illustrating proof of Theorem~\ref{thm:corollary2} for 2d case when
            $A$ matrix is $[a,-b;b,a]$. 
            }\label{fig:illus}
        \end{minipage}
        &
        \begin{minipage}[l]{0.75\textwidth}
            {\small{Let $A$ and $B$ be two adjacent vertices on a regular
            polygon with $2s$ sides and center on origin $O$. Let $C$ be the midpoint of $AB$,
            let $D$ be an arbitrary point on $CB$. 
            Let $\len(AB) = r_0$. 
            If $r$ denotes $\len(OD)$, and $\alpha$ is $\angle{COD}$,
            $r$ as a function of $\alpha$ is given by
            $r(\alpha) = r_0/\cos(\alpha)$. One way to ensure that vector fields point 
            inwards on the boundary is by comparing $dr/d\alpha$. On the trajectories,
            $dr/d\alpha = (dr/dt)/(d\alpha/dt) = ar/b = ar_0/(b\cos(\alpha))$, 
            and on the boundary,
            $dr/d\alpha = r_0\sin(\alpha)/\cos^2(\alpha)$. Vector field points inwards iff
            $|ar_0/(b\cos(\alpha))| \geq |r_0\sin(\alpha)/\cos^2(\alpha)|$, which simplifies to
            $|a/b| \geq \tan(\alpha)$ and maximum value of $\alpha$ is $360/4s=90/s$,
            and we get $|a/b|\geq\tan(90/s)$.}}
        \end{minipage}
    \end{tabular}
\end{figure}

We can generalize Theorem~\ref{thm:corollary} to also synthesize invariants when 
a complex eigenvalue $a + b\iota$ violates $a + |b| < 0$ -- it will involve adding more
faces by increasing $s$.
\begin{theorem}[Generalized sufficient condition for $2s$-LI]\label{thm:corollary2}
The linear system $d{x}/dt = A {x}$,
where $A$ is a $n\times n$ matrix, has a $2s$-LI 
for a value $s$ obtained by adding $\num(\lambda)$ for each distinct eigenvalue
    $\lambda$ of $A$, where:
\begin{itemize}
  \item[(a)] $\num(\lambda) = \mathtt{algm}(\lambda)$ if $\lambda$ is real and $\lambda < 0$
  \item[(b)] $\num(\lambda) = \mathtt{geom}(\lambda)$ if $\lambda$ is real and $\lambda = 0$
  \item[(c)] $\num(\lambda) = k*\mathtt{algm}(\lambda)$ if $\lambda := a + \iota b$ is complex, 
      $a < 0$, and $k > \frac{90}{\tan^{-1}(|a|/|b|)}$
  \item[(d)] $\num(\lambda) = k*\mathtt{geom}(\lambda)$ if $\lambda := a + \iota b$ is complex, 
      $a < 0$, and $k = \frac{90}{\tan^{-1}(|a|/|b|)}$
  \item[(e)] $\num(\lambda) = 0$ otherwise
\end{itemize}
\end{theorem}
\begin{proof}(Sketch)
    The main difference with Theorem~\ref{thm:corollary} is the complex eigenvalue case.
    In this case, we increase the number $s$ of rows in the polyhedral invariant
    depending on the ratio $|a|/|b|$.  The important case to consider is the 2-dimensional
    case when the $A$-matrix is
    $A = [a, b; -b, a]$. Figure~\ref{fig:illus} illustrates this case and shows that
    if $k > \frac{90}{\tan^{-1}(|a|/|b|)}$, then
    any regular polygon with $2k$ sides whose center is the origin will be a 
    $2k$-LI for this 2-dimensional system.  \qed
\end{proof}

We next present a necessary condition for existence of $2n$-LI, and examples that show that the 
gap between the sufficient condition, given in Theorem~\ref{theorem:main}, and the
necessary condition, given below in Theorem~\ref{thm:necessary}, is hard to overcome.
\begin{theorem}[Necessary condition for existence of $2s$-LI]\label{thm:necessary}
If the linear system $d{x}/dt = A {x}$ has a bounded $2s$-LI,
then for every eigenvalue $\lambda$ of $A$, it is the case that either
\begin{itemize}
  \item[(a)] $\lambda$ is real and $\lambda < 0$, or
  \item[(b)] $\lambda = 0$ and $\mathtt{algm}(\lambda) = \mathtt{geom}(\lambda)$, or
  \item[(c)] $\lambda := a + \iota b$ is complex, 
      $a < 0$, and  $|a|/|b| > \tan(90/s)$,
  \item[(d)] $\lambda := a + \iota b$ is complex, 
      $a < 0$, $|a|/|b| = \tan(90/s)$, and $\mathtt{algm}(\lambda) = \mathtt{geom}(\lambda)$,
\end{itemize}
\end{theorem}
\begin{proof}(Sketch)
    Since the $2s$-LI is bounded, the $2s$-LI also corresponds to a weak Lyapunov function. 
    Hence, real eigenvalues have to be non-positive. Moreover, if an eigenvalue is 
    zero, but its $\mathtt{algm}(\lambda)$ is different from $\mathtt{geom}(\lambda)$,
    then we can not get a weak Lyapunov function; for example, consider
    $A = [0, 1; 0, 0]$, which has trajectories escaping any bounded region, and hence
    it has no weak Lyapunov function. 
    For the complex case, note that a $2s$-LI can intersect a $2$-dimensional plane
    in a polygon that can have at most $2s$ faces, and that gives rise to condition~(c) and~(d). \qed
\end{proof}

We now present some examples that show there might not be any simple condition that is both
necessary and sufficient for existence of $2s$-LI in terms of eigenvalues.
\begin{example}
    Consider the 2-dimensional linear system from Example~\ref{example-nonsquarePLF}. 
    This 2d system does not have a $4$-LI, but it has a $6$-LI. 
    These results are predicted by our theorems: specifically,
    the sufficient condition from Theorem~\ref{thm:corollary2} implies
    that this $2d$ system will have a $6$-LI,
    whereas the necessary condition from Theorem~\ref{thm:necessary} implies
    that this $2d$ system will {\em{not}} have a $4$-LI.
    So, our results are strong enough to make perfect predictions for this case.
    
    Next, consider the linear system from Example~\ref{example:main}. It extends 
    the example from Example~\ref{example-nonsquarePLF}
    to 3 dimensions by adding a 3rd dimension
    $x_3$ whose dynamics is given by $\dot{x_3} = -2 x_3$.
    This 3d system continues to have a $6$-LI.
    Note that $90/\tan^{-1}(|a/b|) = 90/\tan^{-1}(1/\sqrt{3}) = 90/30 = 3$.
    If we use the sufficient condition from Theorem~\ref{thm:corollary2} on the $3d$ system,
    we conclude that the $3d$ system would have a $8$-LI, but it does not help us infer
    existence of a $6$-LI.
    The necessary condition from Theorem~\ref{thm:necessary} says that the
    $3d$ system can not have a $4$-LI. Neither theorem says anything conclusive about 
    existence or non-existence of a $6$-LI.  \qed
\end{example}

\begin{example}
    Consider the linear system whose $A$ matrix is given by 
    \begin{equation}
        A = \left( \begin{array}{rrr} -0.5, & b, & 0 
        \\ -b & -0.5, & 0 \\   0, & 0, & -\lambda \end{array} \right)
    \end{equation}
    where $b$ and $\lambda$ are two parameters.
    When $b = \sqrt{3}/2$, then 
    setting $\lambda = 2$ gives us the $3d$ linear system from Example~\ref{example:main}. 
    In fact, when $b = \sqrt{3}/2$, then 
    $\lambda = 2$ is the only choice for $\lambda$ that allows the system to have
    a $6$-LI.  As $b$ decreases toward $1/2$, there are more choices of
    $\lambda$ that imply existence of $6$-LI, but these choices remain in a bounded range
    of values (always upper-bounded by $2$). This simply shows that there is no simple
    property of individual eigenvalues that characterizes existence of $2n$-LI. \qed
\end{example}

\begin{example}
    Consider the $6d$ linear system with $A$ matrix whose element are given by
    \begin{equation}
        a_{ij}  =  \left\{ 
        \begin{array}{rl}  1 & \mbox{ if } i=1,j=2
                        \\ -1 & \mbox{ elsif } i=j \vee j=i+1 \vee i=6,j=1
        \\ 0 & \mbox{ otherwise } \end{array} \right.
    \end{equation}
    This $A$ matrix has eigenvalues $-0.134 \pm 0.5\iota$,
    $-1.866 \pm 0.5\iota$ and $-1\pm \iota$. The value $\tan^{-1}(|a/b|)$ 
    for the eigenvalue $-0.134 + 0.5\iota$
    is roughly 15 degrees, so Theorem~\ref{thm:necessary} tells us that
    this system can not have a $2s$-LI where $s < 90/15=6$. It clearly has a $12$-LI, one 
    such invariant is given by $\max(|x_1|,\ldots,|x_6|) \leq 1$.
    However, our sufficient condition in Theorem~\ref{thm:corollary2} only guarantees existence
    of $2s$-LI where $s$ is at least $6 + 2 + 2 = 10$.

    If we replace $6$ by $8$ in the above definition of $a_{ij}$, we get an $A$ matrix 
    for a $8d$ linear system. This matrix has an eigenvalue $-0.076+0.383\iota$, and
    $\tan^{-1}(0.076/0.383)$ is roughly $11.25$ degrees, and thus Theorem~\ref{thm:necessary}
    implies we need $s$ to be at least $90/11.25 = 8$ to get an $2s$-LI. Indeed there is an
    $16$-LI.  \qed
%
\end{example}

\section{Conclusion}

We presented necessary and sufficient conditions for existence of linear
invariants for linear dynamical systems.  The proof of sufficiency is
constructive and yields a procedure for synthesizing linear invariants
that only needs computation of the eigenvalues and eigenvectors of the
$A$ matrix.  We also presented examples that show the conditions are
tight when applied to specific linear systems.

Our first sufficient condition for existence of $2n$-LI, which is given
in Theorem~\ref{theorem:main}, can be derived from the sufficient condition
for existence of infinity-norm Lyapunov functions presented by Bitsoris and
Kiendl~\cite{Bitsoris91,Kiendl92:TAC}. However, since we are interested in
invariants (and not contractive invariants or Lyapunov functions), we need
to distinguish the cases when algebraic and geometric multiplicities of an
eigenvalue are equal and when not.
Our generalized sufficient condition in Theorem~\ref{thm:corollary2} is novel
and has not been stated before. The same is also true for the
necessary condition in Theorem~\ref{thm:necessary} and the examples 
showing the gap between the necessary and sufficient conditions.

Apart from improving our understanding of linear systems and infinity-norm
(weak) Lyapunov functions, the results can also be used to build verification
tools for piecewise-linear and hybrid systems that just rely on reasoning 
over linear arithmetic. The question of characterizing matrices similar to
some \bmatrix\ (Definition~\ref{def:Bmatrix}) based on its spectral properties
remains open for future work.


%
%
\bibliographystyle{splncs04}
\bibliography{all}
\end{document}